\documentclass[12pt,twoside]{article}
\title{Strongly Gorenstein projective, injective, and flat modules}
\date{}
\usepackage{amsfonts}
\usepackage{amsmath}
\usepackage{amssymb}
\usepackage{latexsym}
\newtheorem{thm}{\bf Theorem}[section]
\newtheorem{cor}[thm]{\bf Corollary}

\newtheorem{prop}[thm]{\bf Proposition}
\newtheorem{defn}[thm]{\bf Definition}
\newtheorem{defns}[thm]{\bf Definitions}
\newtheorem{rem}[thm]{\bf Remark}
\newtheorem{rems}[thm]{\bf Remarks}
\newtheorem{exmp}[thm]{\bf Example}

\catcode`\ç=13
\defç{\c{c}}
\catcode`\é=13
\defé{\'e}
\catcode`\à=13
\defà{\`a}
\catcode`\è=13
\defè{\`e}
\catcode`\â=13
\defâ{\^a}
\catcode`\ù=13
\defù{\`u}
\catcode`\ê=13
\defê{\^e}
\catcode`\î=13
\defî{\^\i}
\catcode`\ô=13
\defô{\^o}
\newcommand{\field}[1]{\mathbb{#1}}

\newcommand{\Z }{\field{Z}}

\def\proof{{\parindent0pt {\bf Proof.\ }}}

\def\FPD{{\rm FPD}}

\def\pd{{\rm pd}}
\def\fd{{\rm fd}}

\def\g{{\rm G\!-\!dim}}

\def\Im{{\rm Im}}

\def\Ker{{\rm Ker}}

\def\sup{{\rm sup}}

\def\Ext{{\rm Ext}}
\def\Tor{{\rm Tor}}
\def\Hom{{\rm Hom}}

\newcommand{\cqfd}
{\hspace{1cm}
\rule{2mm}{2mm}%
\medbreak%
\par%
}
\begin{document}
\thispagestyle{empty}
\maketitle
\vspace*{-2cm}
\begin{center}{\large\bf Driss Bennis and Najib Mahdou}

\bigskip
\small{Department of Mathematics, Faculty of Science and Technology of 
Fez,\\ Box 2202, University S. M.
Ben Abdellah Fez, Morocco, \\  driss.bennis@caramail.com \\ 
mahdou@hotmail.com}
\end{center}

\bigskip\bigskip
\noindent{\large\bf Abstract.} In this paper, we study a
particular case of Gorenstein projective, injective, and flat
modules, which we call, respectively, strongly Gorenstein
projective, injective, and flat modules. These last three classes
of modules give us a new characterization of the first modules,
and confirm that there is an analogy between the notion of
\textquotedblleft Gorenstein projective, injective, and flat
modules\textquotedblright  and the notion of the usual
\textquotedblleft projective, injective, and flat
modules\textquotedblright.

\begin{section}{Introduction}
Throughout this work, $R$ is a commutative ring with identity element,
and all $R$-modules are unital. If $M$ is any $R$-module, we use $\pd_R(M)$ 
and $\fd_R(M)$ to
denote the usual projective and flat dimensions of $M$,  respectively.\\
 It is convenient to use ``local" to refer to (not necessarily
Noetherian) rings with a unique maximal ideal.\bigskip

In 1967-69, Auslander and Bridger \cite{A1,A2} introduced the
G-dimension for finitely generated modules over Noetherian rings;
denoted by $\g(M)$ where $M$ is a finitely generated module. They
proved the inequality  $\g(M)\leq \pd(M)$, with equality
$\g(M)=\pd(M)$ when $\pd(M)$ is finite. We say that G-dimension is
a refinement of
projective dimension.\\
Several decades later, Enochs, Jenda, and Torrecillas \cite{GoIn,
GoInPj, GoPlat} extended the ideas of Auslander and Bridger, and
introduced three homological dimensions, called Gorenstein
projective, injective, and flat dimensions, which all have been
studied extensively by their founders and by Avramov, Christensen,
Foxby, Franklid, Holm, Martsinkovsky, and Xu among others
\cite{tate,CFH,Fox,GdimCM,HH}. They proved that these dimensions
are similar to (and refinements of) the classical homological
dimensions; i.e., projective, injective, and flat dimensions,
respectively.\bigskip

\indent These three Gorenstein dimensions are constructed, via
resolution, in the same way as the usual homological dimensions
with, respectively, Gorenstein projective, injective, and flat
modules, which are defined as follows:
\bigskip

\begin{defns}[\cite{HH}]\label{def-Gmod}\mbox{}
\begin{enumerate}
\item  An $R$-module $M$ is said to be
Gorenstein projective (G-projective for short), if there exists an
exact sequence of  projective modules $$\mathbf{P}=\
\cdots\rightarrow P_1\rightarrow P_0 \rightarrow P^0
         \rightarrow P^1 \rightarrow\cdots$$ such that  $M \cong \Im(P_0
\rightarrow P^0)$ and such that $\Hom_R ( -, Q) $
leaves the sequence $\mathbf{P}$ exact whenever $Q$ is a projective 
module.\\
The exact sequence $\mathbf{P}$ is called a complete projective
resolution.
        \item  The Gorenstein injective
(G-injective for short) modules are defined dually.
        \item
An $R$-module $M$ is said to be
Gorenstein flat (G-flat for short), if there exists an exact
sequence of  flat modules $$\mathbf{F}=\ \cdots\rightarrow
F_1\rightarrow F_0 \rightarrow F^0 \rightarrow F^1
\rightarrow\cdots$$ such that  $M \cong \Im(F_0 \rightarrow
F^0)$ and such that $-\otimes I $
leaves the sequence $\mathbf{F}$ exact whenever $I$ is an injective 
module.\\
The exact sequence $\mathbf{F}$ is called a complete flat
resolution.
\end{enumerate}
\end{defns}

Almost by definition one has the inclusion: \begin{eqnarray*}
\{projective\ modules\}  &\subseteq&  \{G\!-\!projective\ modules\}
\end{eqnarray*}
The main idea of this paper is to introduce and study an
intermediate class of modules called strongly Gorenstein
projective modules (SG-projective for short),
\begin{eqnarray*}
\{projective\ modules\}  &\subseteq&  \{SG\!-\!projective\ modules\}\\
  &\subseteq&  \{G\!-\!projective\ modules\}
\end{eqnarray*}

These modules are defined by considering all modules and
homomorphisms of the complete resolutions of Definitions
\ref{def-Gmod} (1) above are equal (see Definitions
\ref{DefSGproj}). Similarly, we define the strongly Gorenstein
injective, and flat modules (SG-injective, and SG-flat,
respectively, for short) (see Definitions \ref{DefSGproj} and
\ref{DefSGfla}).\bigskip

\indent The simplicity of these modules manifests in the fact that
they are simple characterizations than their  Gorenstein
correspondent modules (see Propositions \ref{caraGfor},
\ref{Gfortf} and \ref{caraGplatfor}  and Remarks \ref{caraGinj}
(2)). Moreover, with these kinds of modules, we are able to give a
nice new characterizations of Gorenstein projective, injective,
and flat modules, similar to the characterization of projective
modules by the free modules, which is the main result of this
paper (see Theorems \ref{caraGproj} and \ref{caraGpla}):\bigskip

\noindent\textbf{Theorem.}\textit{ A module is Gorenstein
projective (resp., injective) if, and only if, it is a direct
summand of a strongly Gorenstein projective (resp.,
injective) module.\\[0.2cm]
Every flat module is  a direct summand of a strongly Gorenstein
flat module.}\bigskip

Over Noetherian rings the Gorenstein projective, injective, and
flat modules were (and still are) excessively studied (please see
\cite{LW}). So, we find that the relation that exists between the
Gorenstein projective and Gorenstein flat modules is (nearly)
similar to the one between the classical  projective and flat
modules (see \cite[Proposition
5.1.4]{LW}\,\footnote{\label{Footnot1}In \cite{LW} Christensen
forgot some details in few results (as \cite[Proposition
5.1.4]{LW}). For the correction, see errata on the Christensen's
homepage:
\begin{center}
http://www.math.unl.edu/$\sim$1christensen3/publications.html\end{center}
} and \cite[Theorem 5.1.11]{LW}). In \cite{HH}, Holm extended
\cite[Proposition 5.1.4]{LW} to coherent rings with finite
finitistic projective dimension. Recall the finitistic projective
dimension of a ring $R$, $\FPD(R)$,    is defined by:
$$\FPD(R)=\sup\{\pd_R(M)|M\;R\!-\!module\;with\;\pd_R(M)<\infty\}$$

\begin{prop}[\cite{HH}, Proposition 3.4] If $R$ is
coherent  with finite finitistic projective dimension, then every
Gorenstein projective $R$-module is  Gorenstein flat.
\end{prop}

Also, \cite[Theorem 5.1.11]{LW} can be extended to coherent rings.
In fact, using Holm's work \cite{HH}, the same proof of
\cite[Theorem 5.1.11]{LW} and \cite[Lemma 5.1.10]{LW} (please see
the footnote p. \pageref{Footnot1}) implies the desired extension,
that is:

\begin{prop}\label{Coh-Gpro-Fla} If $R$ is coherent, then a finitely 
presented $R$-module
is Gorenstein flat if and only if, it is Gorenstein projective.
\end{prop}

In this context, the strongly Gorenstein projective and flat
modules give us more relations. And we prove the two following
results (see Proposition \ref{SGpro-SG-fla-tf} and Corollary
\ref{SGpro-SG-fla-tf2}):\bigskip

\noindent\textbf{Proposition.} \textit{ A module is finitely
generated strongly Gorenstein projective if, and only if, it is
finitely presented strongly Gorenstein flat.}\bigskip

\noindent\textbf{Corollary.} \textit{ If $R$ is integral domain or
local, then a  finitely generated  $R$-module is strongly
Gorenstein flat if, and only if, it is strongly Gorenstein
projective.}\bigskip

The study of finitely generated strongly Gorenstein projective and
flat module allows us to give a new characterization of S-rings.\\
Recall that a ring $R$ is called an S-ring if every finitely
generated flat $R$-module is projective (see \cite{Srin1}). We
have  Proposition \ref{Sring}:\bigskip

\noindent\textbf{Proposition.} \textit{ $R$ is an S-ring if, and
only if, every finitely generated strongly Gorenstein flat
$R$-module is  strongly Gorenstein projective.}\bigskip

Finally, to give credibility to our study, we set some examples to
distinguish the strongly Gorenstein projective, injective, and
flat modules from their correspondent Gorenstein and classical
modules.\bigskip

\end{section}
\begin{section}{Strongly Gorenstein projective and strongly Gorenstein 
injective modules}

\indent In this section we introduce and study the strongly
Gorenstein projective and injective modules which are defined as
follows:

\begin{defn}\label{DefSGproj}
A complete projective resolution of the form  $$ \mathbf{P}=\ 
\cdots\stackrel{f}{\longrightarrow}P\stackrel{f}{\longrightarrow}P\stackrel{f}{\longrightarrow}P
\stackrel{f}{\longrightarrow}\cdots
    $$
   is called strongly complete projective resolution and denoted by
   $(\mathbf{P},f)$.\\[0.2cm]
An  $R$-module $M$ is called strongly Gorenstein projective
(SG-projective for short) if $M\cong \Ker\, f$ for some strongly
complete projective resolution   $(\mathbf{P},f)$.\\[0.2cm]
The strongly Gorenstein injective (SG-injective for short) modules
are defined dually.
\end{defn}

Using the definitions, we immediately get the following results.

\begin{prop}\label{sumSgPro}\begin{enumerate}
    \item If $(P_i)_{i\in I}$ is a family of strongly Gorenstein projective
modules, then $\oplus P_i$ is strongly Gorenstein projective.
    \item If $(I_i)_{i\in I}$ is a family of strongly Gorenstein injective
modules, then $\prod I_i$ is strongly Gorenstein injective.
\end{enumerate}
\end{prop}

\proof Simply note that a sum (resp., product) of strongly
complete projective (resp., injective) resolutions is also a
strongly complete projective (resp., injective) resolution (using
the natural  isomorphisms   in \cite[Theorems 2.4 and 2.6]{Rot}
and \cite[\S 2, N$^\mathrm{o}$ 2, Proposition 1]{Bou}). \cqfd
\bigskip

It is straightforward that the strongly Gorenstein projective
(resp., injective) modules are a particular case of the Gorenstein
projective (resp., injective) modules. And, it is well-known  that
every projective (resp., injective) module is Gorenstein
projective (resp., injective). That is obtained easily by
considering for a projective module $P$
the complete projective resolution
$0\longrightarrow P\stackrel{=}{\longrightarrow}P\longrightarrow0$ 
\cite[Observation 4.2.2]{LW}.\\
\indent Next result shows that the class of all  strongly
Gorenstein projective (resp., injective) modules is between the
class of all projective (resp., injective) modules and the class
of all Gorenstein projective (resp., injective) modules.

\begin{prop}\label{ProSGpro}
Every  projective   (resp., injective) module  is strongly
Gorenstein projective (resp., injective).
\end{prop}

\proof It suffices to prove the Gorenstein projectivity case, and the 
Gorenstein injectivity case is analogous.\\[0.2cm]
Let $P$ be a  projective $R$-module, and  consider the exact
sequence:
$$\begin{array}{ccccccccccc}
    \mathbf{P}= &  & \cdots & \stackrel{f}{\longrightarrow} & P\oplus P & 
\stackrel{f}{\longrightarrow} & P\oplus P&  \stackrel{f}{\longrightarrow} 
&P\oplus P& \stackrel{f}{\longrightarrow}& \cdots \\
         &  &        &                        & (x,y) &\longmapsto & (0,x)  
&   &  &      &
  \end{array}
    $$
We have $0\oplus P=\Ker\, f=\Im \,f\cong P$.\\
Consider a projective module $ Q$, and applying the functor
$\Hom_R(-,Q)$ to the above sequence $\mathbf{P}$, we get the
following commutative diagram
    $$
    \begin{array}{ccccccc}
      \cdots & \longrightarrow & \Hom(P\oplus P, Q) & 
\stackrel{\Hom_R(f,Q)}{\longrightarrow} & \Hom(P\oplus P, Q) 
&\longrightarrow& \cdots \\
        &   &\cong\downarrow &  & \cong\downarrow  &  & \\
      \cdots & \longrightarrow & \Hom(P, Q)\oplus \Hom(P, Q) & 
{\longrightarrow} & \Hom(P, Q)\oplus \Hom(P, Q) & \longrightarrow & \cdots
    \end{array}
    $$
Since the lower sequence in the diagram above is exact, the
proposition follows. \cqfd\bigskip

The strongly Gorenstein projective (resp., injective) modules are
not necessarily projective (resp., injective), as shown by the
following examples. Before that, recall that a ring $R$ is called
quasi-Frobenius (QF-ring for short), if it is Noetherian and
self-injective (i.e., $R$ is an injective $R$-module). For
instance,  if $I$ is a nonzero ideal in a Dedekind domain $R$,
then $R/I$ is quasi-Frobenius \cite[Exercice 9.24]{Rot}. The
following gives a characterization of such rings:

\begin{thm}[\cite{AF}, Theorem 31.9]\label{cariQFring}
The following conditions are equivalent:
\begin{enumerate}
    \item $R$ is quasi-Frobenius;
    \item Every projective $R$-module is injective;
    \item Every injective $R$-module is projective.
\end{enumerate}
\end{thm}

Now we can give the desired examples.

\begin{exmp}\label{ExmSGPro1}
  Consider  the  quasi-Frobenius local ring
$R=k[X]/(X^{2})$ where $k$ is a field, and denote       $\overline{X}$   the 
residue
class in $R$  of $X$.
\begin{enumerate}
    \item  The ideal $(\overline{X})$ is simultaneously  strongly Gorenstein 
projective and strongly Gorenstein
    injective.
    \item But, it is neither  projective nor injective.
\end{enumerate}
  \end{exmp}

\proof 1. With   the homothety $x$ given by
             multiplication by  $\overline{X}$
we have the exact sequence
             $ \mathbf{F}= \, \cdots \longrightarrow R
             \stackrel{x}{\longrightarrow} R
             \stackrel{x}{\longrightarrow} R \longrightarrow\cdots
             $.  Then, $\Ker\, x=\Im\,x=(\overline{X})$.\\
             Since $R$ is quasi-Frobenius, we can see easily from
             Theorem \ref{cariQFring} that $\mathbf{F}$ is
             simultaneously  strongly complete projective and  injective 
resolution. Thus,
               $(\overline{X})$ is both strongly Gorenstein projective and 
injective
               ideal.\\[0.2cm]
\noindent 2.  The ideal $(\overline{X})$  it is not projective,
               since  it is not a free ideal in the local ring $R$ (since
               $\overline{X}^{2}=0$).
               Then, from Theorem \ref{cariQFring} we
               conclude that $\overline{X}$ is also not injective, as
               desired.  \cqfd

\begin{rem}\label{ExmSGPro3}
If we want to construct an example of non finitely generated
strongly Gorenstein projective module,  one can see easily, from
Proposition  \ref{sumSgPro} and using the ideal
$(\overline{X}_{0})$ of the previous example,  that the direct sum
$(\overline{X}_{0})^{(I)}$ for any infinite index set $I$
  is a non finitely generated strongly  Gorenstein projective module.
\end{rem}

Now we give our main result of this paper in which we give a new
characterization of the Gorenstein projective (resp., injective)
modules by the strongly Gorenstein projective (resp., injective)
modules.

\begin{thm}\label{caraGproj}
A module is Gorenstein projective (resp., injective) if, and only
if, it is a direct summand of a strongly Gorenstein projective
(resp., injective) module.
\end{thm}

\proof It suffices to prove the Gorenstein projectivity case, and the 
Gorenstein injectivity case is analogous.\\[0.2cm]
By \cite[Proposition 2.5]{HH}, it remains to prove the direct  
implication.\\[0.2cm]
Let $M$ be a Gorenstein projective module. Then, there exists a
complete projective resolution
$$
\mathbf{P}= \cdots\longrightarrow
P_{1}\stackrel{d^{P}_{1}}{\longrightarrow}
   P_{0}\stackrel{d^{P}_{0}}{\longrightarrow}P_{-1}
   \stackrel{d^{P}_{-1}}{\longrightarrow}
   P_{-2}\longrightarrow\cdots
    $$
    such that $M\cong \Im(d^{P}_{0})$.\\
  For all $ m\in \Z$, denote $\Sigma^{m}P$ the exact sequence
  obtained from $ \mathbf{P}$ by increasing all index by $m$:
$$(\Sigma^{m}P)_{i}=P_{i-m}\quad \
\mathrm{and}\  \quad d_{i}^{\Sigma^{m}P}=d^{P}_{i-m}\qquad\  \mathrm{for}\ 
\mathrm{all}\ i\in \Z.$$
Considering the exact sequence  $$\mathbf{Q}=\oplus(\Sigma^{m}P)=
\cdots\longrightarrow Q=\oplus P_{i} \stackrel{\oplus
d^{P}_{i}}{\longrightarrow}
Q=\oplus P_{i}\stackrel{\oplus
d^{P}_{i}}{\longrightarrow}Q=\oplus P_{i}\longrightarrow\cdots
    $$
Since $\Im(\oplus d_{i})\cong \oplus \Im \,d_{i} $,  $M$ is a direct summand 
of $\Im(\oplus d_{i})$.\\
Moreover, from \cite[Proposition 20.2 (1)]{AF}
$$\Hom(\bigoplus\limits_{m\in \Z}(\Sigma^{m}P),L)\cong \prod\limits_{m\in 
\Z}
    \Hom(\Sigma^{m}P,L)$$ which is an exact sequence for any projective 
module
    $L$.
Thus, $\mathbf{Q}$ is a strongly complete projective resolution.\\
Therefore, $M$
   is a direct summand of the strongly Gorenstein projective module 
$\Im(\oplus
    d_{i})$, as desired.\cqfd

\begin{rem} From \cite[Proposition 2.4]{HH}, we can consider all modules of
the complete projective resolution  in the previous proof are
free, then so are the modules in the constructed strongly complete
projective resolution.
\end{rem}

In the end of this section we give an example of a Gorenstein
projective module which is not strongly Gorenstein projective.
Before that, we give some properties of the strongly Gorenstein
projective modules.\bigskip

The next result gives a simple characterization of the strongly
Gorenstein projective modules.

\begin{prop}\label{caraGfor}
For any module $M$, the following are equivalent:
\begin{enumerate}
  \item $M$ is  strongly Gorenstein projective;
  \item There exists a short exact sequence
$0\rightarrow M\rightarrow P\rightarrow M\rightarrow 0$,
where $P$ is a projective  module,  and  $\Ext(M,Q)=0$ for any   projective 
module $Q$;
  \item There exists a short exact sequence
$0\rightarrow M\rightarrow P\rightarrow M\rightarrow 0$,
where $P$ is a projective  module,  and  $\Ext(M,Q')=0$ for any   module 
$Q'$ with finite projective
dimension;
\item There exists a short exact sequence
$0\rightarrow M\rightarrow P\rightarrow M\rightarrow 0$,
where $P$ is a projective  module;  such that,  for any   projective module 
$Q$, the short sequence   $0\rightarrow \Hom(M,Q)\rightarrow 
\Hom(P,Q)\rightarrow \Hom(M,Q)\rightarrow 0$ is exact;
  \item There exists a short exact sequence
$0\rightarrow M\rightarrow P\rightarrow M\rightarrow 0$,
where $P$ is a projective  module;  such that,  for any module $Q'$ with 
finite projective dimension, the short sequence
$0\rightarrow \Hom(M,Q')\rightarrow \Hom(P,Q')\rightarrow 
\Hom(M,Q')\rightarrow 0$
is exact.
\end{enumerate}
\end{prop}

\proof Using standard arguments, this follows immediately from the
Definition of strongly Gorenstein modules.\cqfd

\begin{rems}\label{caraGinj}\begin{enumerate}
    \item Note that using this characterization of strongly Gorenstein 
projective modules, the Proposition
\ref{ProSGpro} becomes straightforward. Indeed, we have the short
exact sequence  $0\rightarrow P\rightarrow P\oplus P\rightarrow
P\rightarrow
    0$, and  $\Ext(P,Q)=0$ for any module $Q$.
    \item We can also characterize the strongly Gorenstein injective modules 
in a
way similar to the description of strongly Gorenstein projective
modules in Proposition \ref{caraGfor}.
\end{enumerate}
\end{rems}

Recall that a strongly Gorenstein projective module is projective
if, and only if, it has finite projective dimension
\cite[Proposition 2.27]{HH}. In the next result we give similar
result in which the strongly Gorenstein projective modules link
with the flat dimension.

\begin{cor}
A strongly Gorenstein projective module is flat if, and only if,
it has finite flat dimension.
\end{cor}
\proof This is a simple consequence of Proposition
\ref{caraGfor}.\cqfd \bigskip

The following  proposition handles the finitely generated strongly
Gorenstein projective modules. It is well-Known that a finitely
generated projective module is infinitely presented (i.e., it
admits a free resolution $$\cdots \rightarrow F_n\rightarrow
F_{n-1}\rightarrow \cdots \rightarrow F_0\rightarrow M \rightarrow
0$$ such that each $F_i$ is a finitely generated free module).\\
For the Gorenstein projective modules the question is still open.
However, the strongly Gorenstein projective modules give the
following partial affirmative answer, in which we give a
characterization of the finitely generated strongly Gorenstein
projective modules.

\begin{prop}\label{Gfortf} Let $M$ be an $R$-module. The following are 
equivalent:
\begin{enumerate}
\item $M$ is finitely generated strongly Gorenstein projective;
\item There exists a short exact sequence   $0\rightarrow M \rightarrow P
\rightarrow M \rightarrow 0$ where  $P$ is a finitely generated
projective $R$-module, and   $\Ext(M,R)=0$;
\item There exists a short exact sequence   $0\rightarrow M \rightarrow P
\rightarrow M \rightarrow 0$ where  $P$ is a finitely generated
projective $R$-module, and   $\Ext(M,F)=0  $ for all flat
$R$-modules  $F$;
\item  There exists a short exact sequence   $0\rightarrow M \rightarrow P
\rightarrow M \rightarrow 0$ where  $P$ is a finitely generated
projective $R$-module, and $\Ext(M,F')=0  $ for all $R$-modules
$F'$ with finite flat dimension.
\end{enumerate}
\end{prop}

\proof Note that the forth condition is stronger than the first,
this leaves us three implications to prove.\\[0.2cm]
$(1) \Rightarrow (2)$.  This is a simple
consequence of Proposition \ref{caraGfor}.\\[0.2cm]
$(2) \Rightarrow (3)$. Let  $F$ be a flat $R$-module. By Lazard's
Theorem   \cite[\S 1, N$^\mathrm{o}$ 6, Theorem 1]{Bou}, there is
a direct  system  $(L_{i})_{i\in I}$ of finitely generated free
$R$-modules such that $\ \underrightarrow{lim}\;L_{i}\, \cong \,
F$. From \cite[Theorem 2.1.5 (3)]{Glaz}, we have:
\begin{eqnarray*}
\Ext(M,F)&\cong&\Ext(M,\underrightarrow{lim}\;L_{i}) \\
         &\cong&\underrightarrow{lim}\;\Ext(M,L_{i})
\end{eqnarray*}
Now, combining \cite[Theorem 2.1.5 (3)]{Glaz} with \cite[Example
20', page 41]{Rot} shows immediately that
$\Ext(M,L_{i})=0$ for all $i \in I$, as desired.\\[0.2cm]
$(3) \Rightarrow (4)$.  Let $F'$ be an $R$-module such that 
$0<\fd(F')=m<\infty$.\\
  First, we can see easily that (3) implies $\Ext^{n}(M,F)=0$ for all $n>0$, 
and all flat $R$-modules $F$.\\
  Now, pick  a short exact sequence
$0 \rightarrow K \rightarrow L \rightarrow F' \rightarrow 0$
  where  $L$  is a free $R$-module and
   $\fd(K)=m-1$. By  induction $\Ext^{n}(M,L)=\Ext^{n}(M,K)=0$ for all $n> 
0$.  Then, applying the functor $\Hom(M,-)$
   to the  short exact sequence above we obtain the exact sequence
    $$0=\Ext(M,L)\rightarrow \Ext(M,F')\rightarrow
   \Ext^2(M,K)=0$$ Therefore,  $\Ext(M,F')=0$.\cqfd\bigskip

We finish this section by an example of a Gorenstein projective
module which is not strongly Gorenstein projective.

\begin{exmp}\label{ex3}
Consider  the Noetherian local ring
$R=k[[X_{1},X_{2}]]/(X_{1}X_{2})$ with $k$ is a field. Then:
\begin{enumerate}
    \item The two ideals $(\overline{X_{1}})$ and $(\overline{X_{2}})$
are Gorenstein projective, where $\overline{X_{i}}$ is the residue
class in $R$ of $X_{i}$ for $i=1,\, 2$.
    \item $(\overline{X_{1}})$ and $(\overline{X_{2}})$ are not strongly 
Gorenstein projective.
\end{enumerate}
\end{exmp}

\proof 1. This is \cite[Example 4.1.5]{LW}.\\[0.2cm]
\noindent 2. Assume, for example, that the ideal
$(\overline{X_{1}})$
    is strongly
     Gorenstein projective.\\
By  Proposition \ref{Gfortf}, there exists a short exact
sequence
$$  0\rightarrow (\overline{X_{1}})\rightarrow P\rightarrow
(\overline{X_{1}})\rightarrow 0$$ where  $P$ is a finitely
generated projective module. Since $R$ is local, there exists a
positive integer $n$ such that $P\cong R^{n}$. Thus, we can
rewrite the above short exact sequence as follows :
$$(\Phi):\ 0\longrightarrow(\overline{X_{1}})\longrightarrow
R^{n} \longrightarrow (\overline{X_{1}})\longrightarrow 0.$$ On the
other hand, we can see easily that we have the following short
exact sequence
$$
\begin{array}{ccccccccc}
  0&\longrightarrow &(\overline{X_{2}})&\longrightarrow & R & 
\longrightarrow &  (\overline{X_{1}}) & \longrightarrow & 0 \\
   &                &                    &              &x   & \longmapsto   
  & x\overline{X_{1}}   &                    &  \\
\end{array}$$
Thus, by  Schanuel's lemma \cite[Theorem 3.62]{Rot}, we have
$R^{n}\oplus(\overline{X_{2}})\cong
R\oplus(\overline{X_{1}})$.\\
Tensorising  by $k$, the residue field of $R$, we obtain the
following isomorphism of  \linebreak $k$-vector spaces:
$k^{n}\oplus(k\otimes_{R}(\overline{X_{2}}))\cong
k\oplus(k\otimes_{R}(\overline{X_{1}}))$, and we conclude that
$n=1$.  Therefore, the short exact sequence $(\Phi)$ becomes
$$0\longrightarrow(\overline{X_{1}})\stackrel{g}\longrightarrow
R \stackrel{f}\longrightarrow(\overline{X_{1}})\longrightarrow0$$
Now, consider $f(1)=\overline{\alpha}\overline{X_{1}}$ for some
$\alpha \in  k[[X_{1},X_{2}]]$, hence
$\Im\,f=(\overline{\alpha}\overline{X_{1}})=(\overline{X_{1}})$, which 
implies that
  there exist
  $\beta$ and $\delta$ in $k[[X_{1},X_{2}]]$  such that
   $X_{1}=\beta\alpha X_{1}  +\delta
X_{1}X_{2}$,  hence $\alpha \beta=1-\delta X_{2}$ which is
invertible in $k[[X_{1},X_{2}]]$, then so is $\alpha$, and hence
$\overline{\alpha}$ is invertible in $R$. Thus,
$$\Ker\,f=\{x\in
R|xf(1)=x\overline{\alpha}\overline{X_{1}}=\overline{0}\}=\{x\in
R|x\overline{X_{1}}=\overline{0}\}=Ann{\overline{X_{1}}}=(\overline{X_{2}}).$$
Consequently, $(\overline{X_{1}})\cong
\Im\,g=\Ker\,f=(\overline{X_{2}})$.\\
But, this is absurd since  $Ann\,{\overline{X_{1}}}=(\overline{X_{2}})\neq
(\overline{X_{1}})=Ann\,{\overline{X_{2}}}$. \\ Therefore,
$(\overline{X_{1}})$ is not strongly Gorenstein projective.\cqfd
\bigskip

\end{section}

\bigskip
\begin{section}{Strongly Gorenstein flat modules}
\mbox{}\indent In this section we introduce and study the strongly
Gorenstein flat  modules, and further we link them with the
strongly Gorenstein projective modules.
\begin{defn}\label{DefSGfla}
A complete flat resolution of the form  $$ \mathbf{F}=\ 
\cdots\stackrel{f}{\longrightarrow}F
\stackrel{f}{\longrightarrow}F\stackrel{f}{\longrightarrow}F
\stackrel{f}{\longrightarrow}\cdots
    $$
   is called strongly complete flat resolution and denoted by
   $( \mathbf{F},f)$.\\[0.2cm]
An  $R$-module $M$ is called strongly Gorenstein flat (SG-flat for
short) if $M\cong \Ker\, f$ for some strongly complete flat
resolution   $(\mathbf{F},f)$.
\end{defn}

Consequently, the strongly Gorenstein flat modules  is  simple
particular cases of Gorenstein flat modules. The Example \ref{ex4}
gives an example of Gorenstein flat modules which are not strongly
Gorenstein flat.\bigskip

Now, similarly to Proposition \ref{ProSGpro} we prove the following:

\begin{prop}\label{flaSGfla}
Every flat module is strongly Gorenstein flat.
\end{prop}

\begin{exmp}
From   Example \ref{ExmSGPro1}, we can see easily that the ideal
$(\overline{X})$ is also strongly Gorenstein flat, but it is not flat.
\end{exmp}

\begin{prop}\label{sumSgflat}
Every direct sum of strongly Gorenstein flat modules is also
strongly Gorenstein flat.
\end{prop}

\proof Immediate as  the proof of Proposition \ref{sumSgPro} using
the fact that tensorproducts commutes with sums.\cqfd\bigskip

With strongly Gorenstein flat modules we have a simple
characterization of Gorenstein flat modules, that is:

\begin{thm}\label{caraGpla}
If a module is Gorenstein flat, then it is a direct summand of a
strongly Gorenstein flat module.
\end{thm}

\proof  Similar  to the proof of Theorem \ref{caraGproj}.\cqfd
\bigskip

Also, similarly to Proposition \ref{caraGfor}, we have the
following characterization of the strongly Gorenstein flat
modules.

\begin{prop}\label{caraGplatfor}
For any module $M$, the following are equivalent:
\begin{enumerate}
  \item $M$ is  strongly Gorenstein flat;
  \item There exists a short exact sequence
$0\rightarrow M\rightarrow F\rightarrow M\rightarrow 0$,
where $F$ is a flat  module,  and  $\Tor(M,I)=0$ for any injective
module  $I$;
  \item There exists a short exact sequence
$0\rightarrow M\rightarrow F\rightarrow M\rightarrow 0$,
where $F$ is a flat  module,  and  $\Tor(M,I')=0$ for any
module  $I'$ with finite injective
dimension;
\item There exists a short exact sequence
$0\rightarrow M\rightarrow F\rightarrow M\rightarrow 0$,
where $F$ is a flat  module;  such that   the  sequence
$0\rightarrow M\otimes I\rightarrow F\otimes I\rightarrow M\otimes 
I\rightarrow 0$
is exact  for any injective module  $I$;
  \item There exists a short exact sequence
$0\rightarrow M\rightarrow F\rightarrow M\rightarrow 0$,
where $F$ is a flat  module;  such that the   sequence
$0\rightarrow M\otimes I'\rightarrow F\otimes I'\rightarrow M\otimes 
I'\rightarrow 0$ is exact
for any
module  $I'$ with finite injective
dimension.
\end{enumerate}
\end{prop}

Holm \cite[Theorem 3.19]{HH}  proved, over Noetherian rings, that
a Gorenstein flat module is flat if, and only if, it has a finite
flat dimension. Moreover, we can see, from \cite[Proposition
3.11]{HH}, \cite[Theorem 1.2.1]{Glaz}, and the dual of
\cite[Proposition 2.27]{HH} that the same equivalence holds over
coherent rings. But, in general, the question is still open.
However, we can give an other partial affirmative answer
(Corollary 3.8). Before that, we give an affirmative answer in the
case of strongly Gorenstein flat modules.

\begin{prop}\label{Gforpla}
A strongly Gorenstein flat module is flat if, and only if, it has
finite flat dimension.
\end{prop}
\proof Immediate  from Proposition \ref{caraGplatfor}.\cqfd

\begin{cor}
If $R$ has finite weak dimension. Then, an $R$-module is
Gorenstein flat if, and only if, it is flat.
\end{cor}
\proof Simply, combining Theorem \ref{caraGpla} with Proposition
\ref{Gforpla}.\cqfd\bigskip

From Proposition \ref{Coh-Gpro-Fla}, we have that, over coherent
rings, the class of all finitely presented Gorenstein projective
modules and the class of all finitely presented  Gorenstein flat
modules are the same class. In general, the question is still
open. Nevertheless, the strongly Gorenstein modules give the
following partial affirmative answer:

\begin{prop}\label{SGpro-SG-fla-tf}
A module is finitely generated strongly Gorenstein projective if,
and only if, it is finitely presented strongly Gorenstein flat.
\end{prop}

\proof We can prove this similarly to the proof  \cite[Lemma
5.1.10]{LW} using the strongly complete resolutions (please see
the footnote p. \pageref{Footnot1}). Here, we give a proof using
the characterization of  finitely generated strongly Gorenstein
projective modules.\\[0.2cm]
$\Longrightarrow$. Let $M$ be a  finitely generated strongly
Gorenstein projective  module. By Proposition \ref{Gfortf}, there
exists a short exact sequence $0\rightarrow M\rightarrow
P\rightarrow M\rightarrow 0$ where $P$ is a  finitely generated
projective
module, and $\Ext(M,R)=0$.\\
Let $E$ be an injective module. Since $M$ is infinitely presented,
we have, from \cite[Theorem 1.1.8]{Glaz}, the following
isomorphism:
$$\Tor(\Hom(R,E),M)\cong
\Hom(\Ext(M,R),E).$$ Thus, $\Tor(E,M)=0$  (since $\Hom(R,E)\cong
E$ and $\Ext(M,R)=0$). Therefore, $M$ is strongly Gorenstein flat
$R$-module (by Proposition \ref{caraGplatfor}).\\[0.2cm]
$\Longleftarrow$. Now, assume $M$ to be a  finitely presented
strongly Gorenstein flat module. From Proposition
\ref{caraGplatfor}, we deduce that there exists a short exact
sequence $0\rightarrow M\rightarrow P\rightarrow M\rightarrow 0$
where $P$ is a  finitely generated projective module, and
$\Tor(M,E)=0$ for every injective module $E$. If we assume $E$ to
be faithfully injective, the same isomorphism of the direct
implication above implies that $\Ext(M,R)=0$. This means, by
Proposition \ref{Gfortf}, that $M$ is
strongly Gorenstein projective.\cqfd\bigskip

It is well-known that if a flat $R$-module $M$ is
finitely presented, or is finitely generated with $R$ is either local or
integral domain, then $M$ is projective (see \cite[Theorem 3.61
and page 135]{Rot}).\\
Under the same conditions we have the same relation between
strongly  Gorenstein flat modules and strongly  Gorenstein
projective modules, that is Proposition \ref{SGpro-SG-fla-tf} and
the following Corollary:

\begin{cor}\label{SGpro-SG-fla-tf2}
If $R$ is integral domain or local, then a  finitely generated
$R$-module is strongly Gorenstein flat if, and only if, it is
strongly Gorenstein projective.
\end{cor}
\proof Use Proposition \ref{SGpro-SG-fla-tf}.\cqfd\bigskip

Now, we give an example of Gorenstein flat modules which are not
strongly Gorenstein flat.

\begin{exmp}\label{ex4}
Consider the Noetherian local ring $R=k[[X_1,X_2]]/(X_1X_2)$ where
$k$ is a field. Then, the two ideals $(\overline{X_{1}})$ and
$(\overline{X_{2}})$ are Gorenstein flat, where $\overline{X_{i}}$
is the residue class in $R$ of $X_{i}$ for $i=1,2$. But, they are
not strongly Gorenstein flat.
\end{exmp}
\proof Simply apply \cite[Theorem 5.1.11]{LW} and Proposition
\ref{SGpro-SG-fla-tf} to Example \ref{ex3}.\cqfd\bigskip

In studying perfect ring, Bass \cite{bass} proved that a ring $R$
is perfect if, and only if, every flat $R$-module is
projective (see also \cite{AF, Xu} for more details about this ring).\\
Motivated by this result, Sakhajev ask when, more generally, every
finitely generated flat modules is projective (see \cite{srin2}).
In fact, the early study of this question goes back to the 60s,
namely with the considerable works of Vasconcelos \cite{VasSri}
and Endo \cite{Srin4}. However, a first general answer appeared
with Facchini, Herbera, and Sakhajev \cite{Sring}. Recently, an
excessive study of it investigated by Puninski and Rothmaler
\cite{Srin1}, who called the ring which satisfies the question an
S-ring, to honor Sakhajev.\\
Now it is natural to ask: When every finitely generated strongly
Gorenstein  flat module is strongly Gorenstein projective?\\
The answer of this question gives a new characterization of
S-rings, that is:

\begin{prop}\label{Sring}
$R$ is an S-ring if, and only if, every finitely generated
strongly Gorenstein flat $R$-module is  strongly Gorenstein
projective.
\end{prop}
\proof $\Longrightarrow$. Let $M$ be a finitely generated strongly
Gorenstein flat
    $R$-module. Then, by Proposition \ref{caraGplatfor}, there exists a 
short exact sequence
$0\rightarrow M\rightarrow F\rightarrow M\rightarrow 0$
where $F$ is a finitely generated flat $R$-module. By hypothesis
$F$ is projective, and so $M$ is finitely presented. Therefore,
from Proposition \ref{SGpro-SG-fla-tf}, $M$ is strongly Gorenstein 
projective.\\[0.2cm]
$\Longrightarrow$. Now, assume $M$ to be a finitely generated flat
$R$-module. Then, from Proposition \ref{flaSGfla}, $M$ is finitely
generated strongly Gorenstein flat. Hence, it is, by hypothesis,
strongly Gorenstein projective. Thus, from Proposition
\ref{Gfortf},  There exists a short exact sequence $0\rightarrow M
\rightarrow P \rightarrow M \rightarrow 0$ where $P$ is a finitely
generated projective  $R$-module, and $\Ext(M,F)=0 $ for all flat
$R$-modules $F$. Then, $\Ext(M,M)=0 $ (since $M$ is flat), and
then the above short exact sequence split. Therefore, $M$ is
projective as a direct summand of the projective $R$-module $P$,
as desired. \cqfd\bigskip

\noindent {\bf ACKNOWLEDGEMENTS.} The authors would like to express their sincere thanks for 
the referee for his/her helpful suggestions. \\

\end{section}




\bigskip\bigskip

\end{document}